\documentclass[11pt,a4paper,reqno]{article}
\usepackage{caption}
\usepackage{authblk}
\usepackage{amsfonts}
\usepackage{epsfig}
\usepackage{fullpage}
\usepackage{indentfirst}
\usepackage{tikz}
\usepackage{hyperref}
\usepackage{eqnarray,amsthm, amssymb, amsmath,verbatim}
\usepackage{mathrsfs}
\usepackage[margin=1in]{geometry}
\usepackage{enumerate}
\usepackage{graphicx}
\usepackage{hyperref} 
\usepackage{latexsym}
\usetikzlibrary{calc}
\usepackage{framed}
\usepackage{ulem}

\usepackage{enumitem}

\makeatletter
\newtheorem{thm}{Theorem}[section]

\newtheorem{cor}[thm]{Corollary}
\newtheorem*{dfn}{Definition}

\newcommand{\tref}[1]{Theorem \ref{theorem:#1}}

\newcommand{\fref}[1]{Figure \ref{fig:#1}}

\newcommand{\DD}{\mathcal{D}}

\newcommand{\run}{\mathrm{run}}
\newcommand{\runb}{\widetilde{{\mathrm{run}}}}

\newcommand{\rr}{\mathrm{rr}}
\newcommand{\trr}{\widetilde{\mathrm{rr}}}

\newcommand{\signa}{\mathrm{sign}_1}
\newcommand{\signb}{\mathrm{sign}_2}

\newcommand{\initones}{\mathrm{Iones}}

\newcommand{\tDD}{\widetilde{\mathcal{D}}}

\newcommand{\U}{\mathrm{U}}
\newcommand{\area}{\mathrm{area}}
\newcommand{\coarea}{\mathrm{coarea}}

\newcommand{\ret}{\mathrm{ret}}
\newcommand{\lvl}{\mathrm{nc}}
\newcommand{\al}{\alpha}

\newcommand{\la}{\lambda}
\newcommand{\Des}[1]{\mathrm{Des}}

\newcommand{\comp}{\mathrm{comp}}
\newcommand{\rank}{\mathrm{rank}}

\newcommand{\psig}{\Bar{\mathbf{s}}}
\newcommand{\sig}{\mathbf{s}}
\newcommand{\st}{\mathrm{st}}
\newcommand{\cA}{\mathcal{A}}

\newcommand{\thn}[1]{$#1$th}

\newcommand{\colorpath}[3]{
	\draw[line width=2pt, #3] (#1) 
	\foreach \i in {#2} {
		\ifnum\i>0
		-- ++(0,\i)
		\else
		-- ++(-\i,0)
		\fi
	};
}
\usetikzlibrary{calc} % 添加在导言区

\newcommand{\mngrid}[3]{
	\foreach \x in {1,...,#2} {
		\foreach \y in {1,...,#3} {
			\ifnum\numexpr\y*#2 < \numexpr#3*\x+#2\relax
			\ifnum\numexpr\y*#2 > \numexpr#3*(\x-1)\relax
			\path[fill,blue!20!white] ($(#1) + (\x-1,\y-1)$) rectangle +(1,1);
			\fi
			\fi
		}
	}
	\draw[help lines] (#1) grid +(#2,#3);
	\draw[help lines] (#1) -- +(#2,#3);
}

\newcommand{\signseq}[2]{
	\foreach \x/\y [count=\i] in {#2} {
		\node at ($(#1)+(\i-0.5,\x-0.5)$) {$\y$};
	};
}
\newcommand{\signseqb}[2]{
	\foreach \x/\y [count=\i] in {#2} {
		\node[cyan] at ($(#1)+(\i-0.5,\x-0.5)$) {$\y$};
	};
}

\newcommand{\LiLan}{}

\title{Symmetric statistics on rational Dyck paths}
\author[1]{Lilan Dai
	\thanks{\texttt{daililan@mail.nankai.edu.cn}}
}
\author[2]{Shishuo Fu
	\thanks{\texttt{fsshuo@cqu.edu.cn}}
}
\author[3]{Dun Qiu
	\thanks{\texttt{qiudun@nankai.edu.cn}}
}
\affil[1,3]{Center for Combinatorics, LPMC, 
	Nankai University, Tianjin 300071, P. R. China}
\affil[2]{College of Mathematics and Statistics, Chongqing University \& Key Laboratory of Nonlinear Analysis and its Applications (Chongqing University), Ministry of Education, Chongqing 401331, P.R. China}

\date{}

\begin{document}
	\maketitle
	\begin{abstract}
		\noindent
		Rational Dyck paths are the rational generalization of classical Dyck paths. They play an important role in Catalan combinatorics, and have multiple applications in algebra and geometry. In this paper, we introduce two statistics over rational Dyck paths called $\run$ and ratio-run. They both have symmetric joint distributions with the return statistic. We give combinatorial proofs and algebraic proofs of the symmetries, generalizing a result of Li and Lin.\\
        
		
		\noindent\textbf{AMS Classification 2020:} 05A05, 05A15, 05A19
	\end{abstract}

	\section{Introduction}
    
	An $n\times n$ (classical) \textit{Dyck path} is a lattice path from $(0,0)$ to $(n,n)$ consisting of east steps $E$ and north steps $N$ which stays above the diagonal $y=x$. The set of $n\times n$ Dyck paths is denoted by $\DD_n$. It is a well-known result that $|\DD_n|$ is the Catalan number $C_n = \frac{1}{n+1}\binom{2n}{n}$, and Dyck paths are widely studied Catalan objects; see \cite{Sta}.
	
	An $m\times n$ \textit{rational Dyck path} is a lattice path from $(0, 0)$ to $(m,n)$ {staying} {on or above} the line  $y = nx/m$; when $m=kn$, it is called a  $k$\textit{-Dyck path} (or a {\it Fuss--Catalan path}) of size $kn\times n$.
	The set of $m\times n$ rational Dyck paths is denoted by $\DD_{m,n}$. $k$-Dyck paths and rational Dyck paths are generalizations of classical Dyck paths, and have found extensive applications in algebraic combinatorics~\cite{ALW16,Dua2,GMV16} and enumerative geometry~\cite{Sta98} in recent years. For example, $\DD_{kn,n}$ is enumerated by the $k$-Catalan numbers $C_n^{(k)}=\frac{1}{kn+1}\binom{(k+1)n}{n}$, which also count the dominant regions in the $k$-Shi arrangement; see, for example,~\cite{Sta98}.
	
	Two statistics, say $\st_1$ and $\st_2$, defined on the same set of (combinatorial) objects $\cA$ are said to be {\it equidistributed} if the following equation holds%there holds
	\begin{align}\label{equid}
	|\{a\in\cA: \st_1(a)=k\}|=|\{a\in\cA: \st_2(a)=k\}|,
	\end{align}
	for all possible values of $k$. A stronger result called {\it symmetric joint distribution} may be observed for the statistic pair $(\st_1,\st_2)$. It stipulates that for all possible values of $k$ and $\ell$, the following equation holds
	\begin{align}\label{symm equid}
	|\{a\in\cA: \st_1(a)=k,\ \st_2(b)=\ell\}|=|\{a\in\cA: \st_1(a)=\ell,\ \st_2(a)=k\}|.
	\end{align}
	We remark that not only \eqref{symm equid} is aesthetically more desirable than \eqref{equid}, but also that understanding \eqref{symm equid} combinatorially often leads to enlightening bijective/involutive proofs of it and reveals more structural information that $\cA$ bears.

	The joint symmetric distributions of various statistics on classical Dyck paths have already attracted much attention. Using combinatorial bijections, Deutsch \cite{Deu2} proved the symmetry of joint distributions of 
statistics (number of returns, height of the first peak) and (double rises, valleys) on classical 
Dyck paths; Carlsson and Mellit \cite{Car} proved the symmetry of (bounce, area) as a consequence of the shuffle 
theorem; the zeta map $\zeta$ of Haglund \cite{Ha} implies that {\LiLan the pair of statistics (dinv, area) also has a symmetric joint distribution;} Pappe, Paul and Schilling~\cite{Pa} proved the symmetries of (area, depth) and (dinv, ddinv).
    %Using combinatorial bijections, {Deutsch in \cite{Deu2} proved the symmetry of joint distributions of statistics (number of returns, height of the first peak), (doublerises, valleys) on classical Dyck paths};as a consequence of the shuffle theorem, Carlsson and Mellit in \cite{Car} proved the symmetry of the joint distribution of (bounce, area);the zeta map $\zeta$ of Haglund \cite{Ha} implies that the pair of statistics (dinv, area) also have symmetric joint distribution;{Pappe et al.~\cite{Pa} proved the symmetries of the joint distributions of (area, depth) and (dinv, ddinv).}
	On rational Dyck paths, Mellit \cite{Mellit} proved the {symmetric joint distribution for (dinv, area)} as a consequence of the rational shuffle theorem.
	%\redt{Blanco and Petersen \cite{Bla} in 2014 studied the joint distribution of statistics area and rank for Dyck paths.}
	
%	In 2019, motivated by the Pak-Stanley labeling of the chambers of Shi arrangement, Duarte and {Guedes de Oliveira~\cite{Dua2}} introduced the $\run$ statistic for classical Dyck paths. {It counts the number of peaks before the first occurrence of the path pattern $EE$. Let $\ret(P)$ be the number of times that the path $P$ returns to the diagonal. They proved that the statistics $\run$ and $\ret$ are equidistributed on the set of classical Dyck paths.} In 2023, Li and Lin \cite{Li} {strengthened Duarte and Guedes de Oliveira's result and established combinatorially the symmetry of the joint distribution of $(\run, \ret)$. }
	% , which was first observed and established algebraicly by Du~\cite{Du}.

    { 
    Parking functions are Dyck paths such that the north steps are labeled with the numbers in $[n]:=\{1,\ldots,n\}$ and the numbers in each column are increasing from bottom to top. The original definition of parking functions and more related studies on them can be found, for example, in the survey by Yan~\cite[Chapter 13]{Bona}. In the study of hyperplane arrangements, the Pak-Stanley bijection establishes a fundamental correspondence between chambers (described by permutations of $[n]$ decorated with arcs under certain rules) of the Shi arrangement and parking functions; see \cite{Sta96}. The central set of a parking function uniquely determines the arc structure of permutations, thereby enabling the reconstruction of chamber inequalities~\cite{Dua1}. 
    
    In 2019, motivated by the Pak-Stanley labeling of the chambers of Shi arrangement, Duarte and {Guedes de Oliveira~\cite{Dua2}} introduced the \textit{run} statistic for classical Dyck paths, which counts the number of peaks before the first occurrence of the path pattern $EE$. The statistic $\run$ for a parking function is the $\run$ of its underlying Dyck path. They showed that the statistic $\run$ is equidistributed with another statistic called \textit{center}, which implies that chamber inequalities can be equivalently reconstructed through run values. 
    
    For a Dyck path $P$, let the \textit{return} statistic of $P$, denoted by $\ret(P)$, be the number of times that $P$ returns to the diagonal. The return of a parking function is the return of its underlying Dyck path, and a parking function is prime if its return is 1.
    Athanasiadis and Linusson \cite{AL99} gave another bijection between parking functions and regions of the Shi arrangement such that prime parking functions are mapped into bounded regions of the arrangement. 
    This result also suggests that the return statistic of a parking function corresponds to the freeness of the Shi arrangement's regions.
 Furthermore, the return statistic bridges Dyck path geometry and Nakayama algebra homology. As shown in \cite{KlaszMarMar}, it records the repeated projective dimension in the Kupisch series and characterizes isomorphic projective modules.

   %Furthermore, the return statistic is also a bridge between Dyck path geometry and Nakayama algebra homology, see \cite{KlaszMarMar}, that records the repeated projective dimension in the Kupisch series and characterizes isomorphic projective modules.

    In the classical Dyck path case, Du~\cite{Du} proved that the statistics $\run$ and $\ret$ are equidistributed. In 2023, Li and Lin \cite{Li} strengthened Du's result and established combinatorially the symmetry of the joint distribution of $(\run, \ret)$. 

    %In the classical setting, there exists a one-to-many correspondence between Dyck paths and parking functions, where the number of returns in a Dyck path corresponds to the number of prime factors in the parking function's decomposition. A particularly profound combinatorial-geometric connection emerges when examining prime parking functions: these correspond bijectively to the bounded regions of the Shi arrangement \cite{AL99}. Remarkably, parking functions with return statistic equal to 1 are precisely the prime parking functions that parametrize these bounded regions. This observation suggests that the return statistic functions as a fundamental combinatorial invariant, effectively capturing the ``degrees of freedom" and intrinsic structural complexity of the Shi arrangement's bounded regions. 

    Considering the $k$-Dyck path case, both the return and run statistics have natural generalizations to this more general setting. From the hyperplane arrangement aspect, Duarte et al. \cite{Dua3} and Fu et al. \cite{Fu} independently established a bijection between $k$-Dyck paths and $k$-Catalan arrangements.
    These results enable people to explore the properties of the $k$-Catalan arrangements from the statistics defined over $k$-Dyck paths.

    The rational Dyck path case generalizes the $k$-Dyck path case further. In this paper, we define two statistics, {run} and ratio-run (denoted by $\runb$), for the rational case when $m\geq n$, generalizing the run statistic for the classical case. As it turns out, both newly introduced statistics have symmetric distribution with the return statistic over rational Dyck paths. In particular, we have composition type preserving involutions $\Phi$ and $\Psi$ such that
$$
         \run(\Phi(P))=\ret(P),\quad \ret(\Phi(P))=\run(P),
$$
and
$$
         \runb(\Psi(P))=\ret(P),\quad \ret(\Psi(P))=\runb(P).
$$ 
    The above results are summarized as \tref{run} and \tref{runb1} in Section 3. 
    Our results in the rational case suggest that statistics $\run$, $\runb$, and $\ret$ may reveal important properties of other combinatorial and geometric structures, potentially generalizing the results in the Shi case and establishing new connections between lattice paths and hyperplane arrangements.

     %While classical Dyck paths already connect to Nakayama algebras via the Kupisch series, rational Dyck paths encode even richer structures, suggesting that their returns could unlock new perspectives on higher-dimensional or non-commutative generalizations of these algebras. Investigating these returns more deeply may lead to fresh invariants or simplified homological computations in Nakayama algebra theory.
    
	%\begin{thm}\label{theorem:run}
	%	Given positive integers $m,n$ such that $m\geq n$, there is a composition type preserving involution $\Phi$ on the set $\DD_{m,n}$, sending statistics $(\run, \ret)$ to $(\ret, \run)$.
	%\end{thm}
	%\begin{thm}\label{theorem:runb}
	%	Given positive integers $m,n$ such that $m\geq n$, there is a composition type preserving involution $\Psi$ on the set $\DD_{m,n}$, sending statistics $(\runb, \ret)$ to $(\ret, \runb)$.
	%\end{thm}
	
	The organization of this paper is as follows. In Section~\ref{preliminary}, we give necessary backgrounds on rational Dyck paths {and the definitions of the aforementioned statistics.} In Section~\ref{involution},
	we give bijective proofs of the symmetries of the joint distributions $(\run, \ret)$ and $(\runb, \ret)$.
	%given $m>n$, we construct the composition type preserving involutions on the set $\DD_{m,n}$, sending statistics $(\run, \ret)$ to $(\ret, \run)$. And we further provide symmetric joint distribution identity and so on.
	{We establish in Section~\ref{gen} several generating function identities} involving statistics $\ret$ and $\run$ in the classical case and the Fuss--Catalan case. Finally, in Section~\ref{open}, we propose several open problems for future research.

	\section{Preliminaries}\label{preliminary}
	
	Let $P\in\DD_{m,n}$ be an $m\times n$ Dyck path. For the various statistics that we introduce next, the {cells} all refer to those contained in the $m\times n$ rectangle. The \textit{area} of $P$ is the number of cells between the path and the diagonal, denoted by $\area(P)$; the \textit{coarea} of $P$ is the number of cells above the path, denoted by $\coarea(P)$. The \textit{coarea sequence} of $P$ is a sequence of $n$ non-negative integers $\U(P)=(u_1,u_2,... ,u_n)$, where $u_i$ is the number of cells to the left of \thn{i} north step from bottom to top. For example, the Dyck path in \fref{path} has the coarea sequence $(0,0,1,1,3)$. For convenience, when the size of a Dyck path is given, we can use the coarea sequence $(u_1,u_2,...,u_n)$ to represent the path.

	\begin{figure}[ht!]
		\centering
		\begin{tikzpicture}[scale=0.4]
			\mngrid{0,0}{7}{5}
			\colorpath{0,0}{2,-1,2,-2,1,-4}{}
		\end{tikzpicture}
		\caption{A $7\times 5$ Dyck path indexed by coarea sequence $(0,0,1,1,3)$.}
		\label{fig:path}
	\end{figure}
	
	Let the {\it diagonal cells} of a path be the cells that the main diagonal intersects (marked blue in \fref{path}), we are ready to define the statistics return, run and $\runb$ for rational Dyck paths.
	
	\begin{dfn}[return]
		For a path $P\in\DD_{m,n}$, the return of $P$ is the number of {\LiLan $N$ steps that are on the boundary of diagonal cells}, denoted by $\ret(P)$.
	\end{dfn}
	For classical $n\times n$ Dyck paths, the rational return statistic reduces to the classical return statistic.
	
	For any lattice point $(i,j)$ in $m\times n$ lattice, we define the \textit{rank} of $(i,j)$ by $\rank(i,j) {\LiLan:=} j\cdot m - i\cdot n$. {\LiLan It follows that the rank of every point on an $m\times n$ Dyck path is non-negative.} In the graph of an $m\times n$ Dyck path, it tells the distance between the point and the main diagonal.
	Notice that for any $m\times n$ Dyck path with $m\ge n$, the $N$ steps corresponding to returns are in bijection with lattice points on the path (except for the point $(m,n)$) with rank smaller than $n$. We shall call all such lattice points \textit{return positions}.
	
	\begin{figure}[ht!]
		\centering
		\begin{tikzpicture}[scale=0.4]
			\mngrid{0,0}{15}{7}
			\colorpath{0,0}{1,-1,1,-1,2,-1,1,-7,2,-5}{}
			\draw[fill,red] (10,5) circle (5pt);
			\draw[fill,red] (0,0) circle (5pt);
		\end{tikzpicture}
		\caption{A $15\times 7$ Dyck path.}
		\label{fig:fuss}
	\end{figure}
	
	We have two generalizations, run and ratio-run, of the classical run statistic to the rational case, denoted by run and $\runb$ respectively. 
	\begin{dfn}[run]
		For a path $P\in\DD_{m,n}$ with $\U(P)=(u_1,u_2,...,u_n)$, the run statistic is 
		$$\run(P)=\min\{i\in[n]:i\notin\{u_1,u_2,...,u_n\}\}.$$
	\end{dfn}
	In other words, the rational run statistic is the number of peaks before the first occurrence of $EE$ path pattern, the same as the definition of classical run statistic.  We see that these peaks appear in column $1,\ldots,\run(P)$ of $P$. 
Without causing confusion, we also call the steps $N^aE$ that correspond to each of these peaks a \textit{run}. In contrast, the second generalization of $\run$ is defined as follows.
	\begin{dfn}[ratio-run]
		For a path $P\in\DD_{m,n}$ with $\U(P)=(u_1,u_2,...,u_n)$, let $k=\lfloor \frac{m}{n} \rfloor$, then the ratio-run statistic is 
		$$\runb(P)=\min\{i\in[n]:ki\notin\{u_1,u_2,...,u_n\}\}.$$
	\end{dfn}
	Notice that $k$ in the definition is the ratio of $m$ and $n$. 
    The ratio-run statistic is the number $r$ of peaks in lattice columns $1,k+1,2k+1,3k+1,\ldots$ before the first occurrence of $EE$ path pattern in columns $rk$ and $rk+1$. We call the collection of steps confined in columns $ak+1,\ldots,(a+1)k$ a \textit{ratio-run} for each $0\leq a<\runb(P)$. {\LiLan Thus,} all ratio-runs span the columns $1,\ldots,\runb(P)\cdot k$ of $P$.
    % We see that the ratio-runs appear in column $1,\ldots,\runb(P)\cdot k$ of $P$. We call the collection of steps in columns $ak+1,\ldots,(a+1)k$ a \textit{ratio-run} for each $0\leq a<\runb(P)$. 
    
	%All the three definitions are natural generalization of the definition of run and return for classical Dyck paths. Particularly, the first run statistic is  exactly the number of peaks before the first occurrence of EE pattern, same as the definition of run for the classical case.
	
	For example, the $15\times 7$ Dyck path $P$ with coarea sequence $(0,1,2,2,3,10,10)$ in \fref{fuss} has $\ret(P)=2$ and $\run(P)=4$. Return positions are marked in red. The four runs are $NE$, $NE$, $NNE$ and $NE$.
	Since $k=\lfloor \frac{15}{7} \rfloor = 2$, the smallest number $i$ such that $2i$ is not in $\{0,1,2,2,3,10,10\}$ is
	$i=2$, therefore $\runb(P)=2$. The two ratio-runs are $NENE$ and $NNENE$.

	Given an integer $n$, {\LiLan we say that a sequence of positive integers $\al=(\al_1,\ldots,\al_{\ell})$ is a (strong) \textit{composition} of $n$ if $\al_1+\cdots+\al_{\ell}=n$, denoted by $\al\vDash n$, where $|\al|$ denotes the sum of parts}. A \textit{partition} $\la$ of $n$ is a composition such that the parts are ordered non-increasingly, denoted by $\la\vdash n$. The \textit{composition type} $\comp(P) \vDash n$ of a Dyck path $P \in \DD_{m,n}$ is the composition formed by lengths of maximal consecutive north steps of the path, and the \textit{partition type} $\la(P) \vdash n$ is the underlying partition of $\comp(P)$. 
	For example, the Dyck path in \fref{fuss} has composition type $(1,1,2,1,2)$ and  partition type $(2,2,1,1,1)$, respectively. If $\comp(P)=(\al_1,\ldots,\al_\ell)$, we call the subpath $N^{\al_i}E$ the \textit{$i$th vertical component} of $P$; if $\run(P)=k$, we also call the first $k$ vertical components the \textit{run components}. We call the vertical segments corresponding to the return positions the \textit{return components}. We write
	$$\DD_{m,\al}:=\{P\in \DD_{m,|\al|}:\comp(P) = \al\}.$$

	\section{Generalized results on rational Dyck paths}\label{involution}
	
	In this section, we shall state two main theorems about symmetric distributions of statistics and give bijective proofs of them. The bijections in the proofs are involutions of the set $\DD_{m,n}$, exchanging relevant statistics. For our convenience, we use steps (reading from $(0,0)$ to $(m,n)$) to represent a path. For example, the path in \fref{fuss} is $NENENNENEEEEEEENNEEEEE$.
	
	\subsection{\texorpdfstring{The symmetry of $(\run,\ret)$ distribution on $\DD_{m,n}$}{The symmetry of (run,ret) distribution on Dmn}}\label{run-ret}
	{ 
    For the run statistic, we have the symmetry of $(\run,\ret)$ distribution on $\DD_{m,n}$ summarized as the following theorem.
    
    \begin{thm}\label{theorem:run}
    		Given positive integers $m>n$, there is a composition type preserving involution
		\begin{eqnarray*}
			\Phi:& \DD_{m,n} &\longrightarrow \quad\DD_{m,n},\\
			&P &\longmapsto  \quad\Phi(P),
		\end{eqnarray*}
		such that $\Phi(\Phi(P))=P$, $\run(\Phi(P))=\ret(P)$ and $\ret(\Phi(P))=\run(P)$.
		%Given positive integers $m,n$ such that $m\geq n$, there is a composition type preserving involution $\Phi$ on the set $\DD_{m,n}$, sending statistics $(\run, \ret)$ to $(\ret, \run)$.
	\end{thm}}
	
	In order to prove \tref{run}, we shall introduce a few definitions. 
	
	Since there is a classical bijection from $\DD_{n,n}$ to $\DD_{n+1,n}$ sending a path $P$ to $PE$ keeping the statistics run and return unchanged, we shall ignore the case when $m=n$ and only consider when $m>n$.
	
	For any $m\times n$-rational Dyck path such that $m> n$, all runs  occur in the first few columns, while returns occur before {any north step $(i,j)\to (i,j+1)$ with $\rank(i,j)<n$.} 
	
	Given a path $P\in\DD_{m,n}$, suppose that the composition type $\comp(P)=\al=(\al_1,\ldots,\al_\ell)$ and $\run(P)=k$. Let $M_i = N^{\al_i} E$ be the $i$th vertical component of $P$, then 
	$$P=M_1M_2\cdots M_k E E^{\beta_0} M_{k+1} E^{\beta_1} M_{k+2} E^{\beta_2}\cdots M_{\ell} E^{\beta_{\ell-k}}$$
	for some nonnegative integers $\beta_0,\beta_1,\ldots,\beta_{\ell-k}$ such that $\beta_0+\beta_1+\cdots+\beta_{\ell-k} = m-\ell-1$.
	
	In this decomposition, each $M_i$ has one $E$ step and at least one $N$ step. The ending point of $M_i$ has a larger rank than the starting point of $M_i$ since $m>n$. Consequently, there is a unique integer $s$ satisfying $0\leq s\leq k-1$ such that the endpoints
    {\LiLan of each of the first $s$ components $M_i$ has a rank smaller than $n$, while the endpoint of $M_{s+1}$ has a rank greater than or equal to $n$}. There cannot be an $E$ step after the first $s$ $M_i$'s because the rank of the path should be above 0, and they are forced to contribute to both run and return. For $M_{s+1}$, since it is still a run and its start point is a return position, it is also forced to contribute to both run and return.
    We shall let $s+1$ { be} the {\textit{run-return} statistic (or \textit{rr} statistic)} of $P$, representing the number of initial vertical components contributing to both $\run(P)$ and $\ret(P)$, denoted by $\rr(P)$. {\LiLan Therefore, $\rr(P)-1$ is the number of run components that end at a return position.} This also implies the inequality $\rr(P)\le \min\{\ret(P),\run(P)\}$ that holds for every $P\in\DD_{m,n}$.
	
	%Now \tref{run} is equivalent to the following theorem.
	%\begin{thm}\label{theorem:run1}
	%	Given positive integers $m>n$, there is a composition type preserving involution
	%	\begin{eqnarray*}
	%		\Phi:& \DD_{m,n} &\longrightarrow \quad\DD_{m,n},\\
	%		&P &\longmapsto  \quad\Phi(P),
	%	\end{eqnarray*}
	%	such that $\Phi(\Phi(P))=P$, $\run(\Phi(P))=\ret(P)$ and $\ret(\Phi(P))=\run(P)$.
	%\end{thm}
	% We let $\rr(P)$ denote the biggest integer $i$ such that the first $i$ vertical components all contribute to both run and ret statistics, where $\rr(P)-1$ is the number of run components stick close to the diagonal cells.
	
	Let $P$ be a path in $\DD_{m,n}$, it is not hard to see that the rr statistic depends only on the composition type $\comp(P)=\al\vDash n$ and the number $m$.
	Let $\initones(\al)$ denote the length of the maximal initial consecutive singleton segments in the composition $\al$. Then 
	\begin{equation}\label{rr0}
		\rr(P) = \min\left\{\initones(\comp(P))+1,\left\lceil\frac{n}{m-n}\right\rceil\right\}.
	\end{equation}
	In cases where no confusion arises, we denote the statistic $\rr(P)$ by $\rr(\al,m)$, i.e.,\
	\begin{equation}\label{rr}
		\rr(\al,m) = \min\left\{\initones(\al)+1,\left\lceil\frac{|\al|}{m-|\al|}\right\rceil\right\}.
	\end{equation}

    Let 
	$$\DD_{m,n}^{(k,\ell)}:=\{P\in \DD_{m,n}:\run(P)=k, \ret(P)=\ell\},$$
	and
	$$\DD_{m,\al}^{(k,\ell)}:=\{P\in \DD_{m,\al}:\run(P)=k,\ \ret(P)=\ell\}.$$
	In order to prove the main theorem, we shall give a combinatorial proof of the following stronger result. 

	\begin{thm}\label{theorem:run2}
		Let $m>n$ be positive integers and let $\al = (\al_1,\ldots,\al_\ell)$ be a composition of $n$. If $r=\rr(\al,m)$, then
		\begin{equation}\label{runthm}
			\left|\DD_{m,\al}^{(r+a,r+b-a)}\right| = \left|\DD_{m,\al}^{(r+b,r)}\right|
		\end{equation}
		for any non-negative integers $a\leq b$.
	\end{thm}
	\begin{proof}
		We shall prove equation (\ref{runthm}) by giving a series of bijections $\phi$ from $\DD_{m,\al}^{(r+b,r)}$ to $\DD_{m,\al}^{(r+a,r+b-a)}$ for any $a\leq b$.
        
		We introduce the notion of the signature of a lattice path, which is crucial in our proof. For any lattice path $P$, we say that an east step $E$ of $P$ has signature $1$ if it follows a north step, and it has signature $0$ otherwise. The \textit{signature} of $P$, denoted as $\sig(P)$, is the signature sequence of all its east steps when we scan $P$ from its start (southwest) to its end (northeast).
        
        {\LiLan If $P$ is a Dyck path, we let the \textit{plain signature} of $P$, denoted as $\psig(P)$, be the signature sequence of all its east steps
        that are 
        \begin{itemize}
			\item \textit{not} contained in any run or return component, and
			\item \textit{not} following a run component.
		\end{itemize}
        }
       {\LiLan \noindent For example, the two paths in \fref{phieg} have the same plain signature that $\psig(P)=\psig(\phi(P))=0001100000$, while their full signatures are different: $\sig(P)=111100001100000$ and $\sig(\phi(P))=110010011000100$.}
If we know the quadruple (run, return, {\LiLan plain} signature, composition type) of a valid Dyck path $P$, then we know all the information of $P$ since we can reconstruct $P$ by the quadruple. Notice that a random quadruple may not be related to a Dyck path. We will classify all valid {\LiLan plain} signatures for a given triple (run, return, composition type) in Remark 1 at the end of this subsection.

		{\LiLan Fix a path $P\in \DD_{m,\alpha}^{(r+b,r)}$, and write its $i$th vertical component as $M_i = N^{\alpha_i}E$. Since $P$ has exactly $r$ returns, which is minimal for the given composition type $\al$, it never hits any return position after its $r$th run component.} Recall that there is a forced $E$ step after the run components of $P$.
        {\LiLan If we let $\mathcal{T}(m,\alpha,r+b)$ be the set of all paths $T$ from $(r+b+1,\al_1+\cdots+\al_{r+b})$ to $(m,n)$ staying weakly above the main diagonal with no north step bordering a diagonal cell and having composition type $(\al_{r+b+1},\ldots,\al_\ell)$, then we have the decomposition }
        \begin{equation}\label{eq:decomp-P}
        P=M_1M_2\cdots M_{r+b}ET,
        \end{equation}
        {\LiLan where the suffix  path $T\in\mathcal{T}(m,\alpha,r+b)$, creating no new return positions. 
        The set $\mathcal{T}(m,\alpha,r+b)$ is the collection of all such valid suffix paths $T$ of $P\in\DD_{m,\alpha}^{(r+b,r)}$. Define 
        $$
        \sig(m,\alpha,r+b) := \{\sig(T):T\in\mathcal{T}(m,\alpha,r+b)\},
        $$
        and observe that $\psig(P)=\sig(T)$ in the decomposition (\ref{eq:decomp-P}), we conclude that $\sig(m,\alpha,r+b)$ is the set of valid plain signatures of paths in $\DD_{m,\alpha}^{(r+b,r)}$, each element of which is a $01$-sequence of length $m-r-b-1$, with $\ell-r-b$ $1$'s and $m-\ell-1$ $0$'s.}
        
        By following the ``hit and lift'' algorithm that is described and framed below, we derive as output a path $X$ and simply take it to be the image $\phi(P):=X$.

		The main idea is that, thanks to the {\LiLan composition type preserving} restriction that $\comp(P)=\comp(\phi(P))$, for each occurrence of consecutive $EN$ (we shall refer to this intuitively as a ``lift'') in the image path $\phi(P)$, we know precisely how many consecutive copies of $N$ need to be swept before we meet the next $E$ step. Then it is just a matter of knowing where these $EN$ pairs (lifts) take place. The details are summarized as the following algorithm written in a pseudocode fashion. 

		\begin{framed}
		\begin{center}
		 {\bf Hit and Lift Algorithm}
		\end{center}
		\begin{itemize}[leftmargin=5em]
		 {\LiLan \item[{\bf Begin}] Input $m$, $\al=(\al_1,\ldots,\al_\ell)$, $a$, $b$, $P\in D^{(r+b,r)}_{m,\alpha}$, where $r=\rr(\alpha,m)$.
        \item[{\bf Setup}] Decompose $P$ as $P = M_1 M_2\cdots M_{r+b} ET$; \\compute $\psig(P)=\sig(T)=s_1\cdots s_{m-r-b-1}$; \\set $X=M_1M_2\cdots M_{r+a}E$; 
        \\set $i=1$ and $j=r+a+1$.}

		  \item[{\bf Sweep}] While $i\le m-r-b-1$, do:\\
		  			   \phantom{place} if $s_i=1$:\\
                      \phantom{place else} goto ``Lift''; \\
		  			   \phantom{place} else if $XE$ hits a return position:\\
                      \phantom{place else else} set $X=XE$ and goto ``Lift'' (causing an extra Lift);\\
		  			   \phantom{place else }else:\\
                      \phantom{place else else} set $X=XE$, $i=i+1$.
		  \item[{\bf Lift}] While $j\le \ell$, set $X=XN^{\al_{j}}E$, $j=j+1$, $i=i+1$, and go back to ``Sweep''.
		  \item[{\bf End}] Output $X$.
		\end{itemize}
		\end{framed}

		\fref{phieg} shows an example with $m=15$, $\al=(1,2,3,2,1,2)\vDash 11$, and $r=\rr(\al,m)=2$. Set $b=2$ and $a=0$. Then for the path $P\in\DD_{15,(1,2,3,2,1,2)}^{(2+2,2)}$ with {\LiLan plain} signature $0001100000$,  we obtain by the algorithm its image $\phi(P)=X\in\DD_{15,(1,2,3,2,1,2)}^{(2+0,2+2)}$ with the same  {\LiLan plain} signature $0001100000$. The two green circled $1$'s are the results of extra Lifts when hitting the return positions, and they are not counted in {\LiLan plain} signature since they are contained in some return components.
		
		\begin{figure}[ht!]
			\centering
			\begin{tikzpicture}[scale=0.38]
				\mngrid{0,0}{15}{11}
				\colorpath{0,0}{1,-1,2,-1,3,-1,2,-1}{black}
				\draw[red,dotted,line width=2pt] (4,8) -- (5,8);
				\colorpath{5,8}{-3,1,-1,2,-6}{cyan}
				
				\mngrid{19,0}{15}{11}
				\colorpath{19,0}{1,-1,2,-1}{black}
				\colorpath{21,3}{3,-1,2,-1}{black!20}
				\colorpath{19+4,8}{-1}{red!20,dotted}
				\colorpath{24,8}{-3,1,-1,2,-6}{blue!20}
                \signseqb{0,0}{2/1,4/1,7/1,9/1,9/0}
				\signseq{5,8}{1/0,1/0,1/0,2/1,4/1,4/0,4/0,4/0,4/0,4/0}
				
				\node at (17,5.5) {$\longrightarrow$};
				
				\colorpath{19+2,3}{-1}{red,dotted}
				\colorpath{19+3,3}{-1,3,-3,2,-1,1,-4,2,-3}{magenta}
				%		\draw[red,dotted,line width=2pt] (19+4,8) -- +(1,0);
				\draw[green,line width=1pt] (22,3) -- +(1,3);
				\draw[green,line width=1pt] (23,3) -- +(1,3);
				\draw[green,line width=1pt] (30,9) -- +(1,2);
				\draw[green,line width=1pt] (31,9) -- +(1,2);
                \signseqb{19,1}{1/1,3/1,3/0}
				\signseq{22,3}{1/0,4/ ,4/0,4/0,6/1,7/1,7/0,7/0,7/0,9/ ,9/0,9/0}
                    \signseqb{22,3}{1/ ,4/1,4/ ,4/ ,6/ ,7/ ,7/ ,7/ ,7/ ,9/1}
				\draw[green,line width=1pt] (23.5,6.5) circle [radius=0.45];
				\draw[green,line width=1pt] (31.5,11.5) circle [radius=0.45];
			\end{tikzpicture}
			\caption{$P\in\DD_{15,(1,2,3,2,1,2)}^{(4,2)}$ and $\phi(P)\in\DD_{15,(1,2,3,2,1,2)}^{(2,4)}$.}
			\label{fig:phieg}
		\end{figure}
		
		It remains to show that for any pair of non-negative integers $a\leq b$, the map $\phi:\DD_{m,\al}^{(r+b,r)}\rightarrow\DD_{m,\al}^{(r+a,r+b-a)}$ as defined above is 
		\begin{itemize}
			\item well-defined, i.e., $\phi(P)\in\DD_{m,\al}$ and $\run(\phi(P))=r+a$, $\ret(\phi(P))=r+b-a$ for any $P\in \DD_{m,\al}^{(r+b,r)}$;
			\item bijective.
		\end{itemize}
		To show well-definedness, we first explain why the algorithm ensures that the growing path $X$ never crosses the main diagonal $y=nx/m$. Since the slope of the main diagonal $n/m<1$, the vertical sections of the diagonal cells (the blue ribbon strip) are either $1\times 1$ squares or $1\times 2$ dominoes. This means that whenever the path hits a return position, continuing with one $N$ step is enough to ``save it'' from crossing the main diagonal. This is the case since the ``Lift'' step in the algorithm is triggered as soon as the path $X$ hits a return position{\LiLan, thereby inserting} an extra Lift. 

        We shall also claim that the algorithm will produce exactly $b-a$ extra Lifts and end at the lattice point $(m,n)$. 
        {\LiLan Any Dyck path $X$ in $\DD_{m,\al}^{(r+a,r+b-a)}\subseteq\DD_{m,\al}$ has $\ell$ vertical segments since $\al$ has $\ell$ parts, and $m-\ell$ $E$ steps that are not contained in any vertical segments. Recall that $M_i=N^{\al_i}E$ is the $i$th vertical segment with $\al_i\geq 1$ $N$ steps. Since $m>n$, appending a vertical segment $M_i$ during the construction of $X$ increases the rank by $(m\cdot \al_i - n)>0$, while appending an $E$ step decreases the rank by $|-n|>0$.

        It is easier to see that we cannot have more than $b-a$ extra Lifts. When the $(b-a)$th extra Lift is triggered, suppose that we have used $x$ $1$'s and $y$ $0$'s in the plain signature, then we have $\ell-r-b-x$ $1$'s and $m-\ell-1-y$ $0$'s in the plain signature that are not used. We have created $(r+a)+(b-a)+x = r+b+x$ vertical segments, and the rank of the current endpoint is 
        \[
        m\cdot\left(\sum_{i=1}^{r+b+x}\al_i\right)-n\cdot (r+b+x+y+1).\] It is the same with the rank of the origin path $P$ after reading the $(x+y)$th signature position, since the two paths use exactly the same collection of vertical segments and east steps. As a consequence, the remaining part of path $X$ can be constructed directly from the plain signature sequence, will not generate any new extra Lifts, and is identical to the remaining part of path $P$. For example in \fref{phieg}, the second(last) extra Lift (the right circled 1) is triggered on column 13 in the right picture, thus the steps of $P$ and $\phi(P)$ are identical after that Lift on columns 14 and 15.

        We still need to explain that the algorithm cannot terminate with only $i$ extra lifts, where $i<b-a$.
        Suppose, to the contrary, that this occurs and the $(i+1)$-st extra lift is not triggered.
        The algorithm ends after all signatures are used, and the output path $X$ has 
        $(r+a)+i+(\ell-r-b)<(r+a)+(b-a)+(\ell-r-b)=\ell$ vertical segments, while it has the same number of east steps that do not belong to any vertical segment.
        Since every vertical segment contributes positively to the rank,
        $X$ must reach a point whose rank is strictly less than $0$ (the rank of the endpoint of $P$),
        and it drops below the main diagonal, a contradiction. }So $\phi(P)\in\DD_{m,\al}$ and the well-definedness is proved. It also follows immediately that the {\LiLan plain} signature is unchanged.

		Next, we check the values of the statistics. $\run(\phi(P))=r+a$ is evident since $\phi(P)$ has the prefix $M_1M_2\cdots M_{r+a}E$. For the return statistic, note that $\phi(P)$ hits the return positions $r$ times with its run components, and $b-a$ times afterwards when extra Lifts occur, corresponding to the instances when the algorithm runs into the ``else if'' case. So in total $\phi(P)$ hits the return positions precisely $r+b-a$ times. 

		% notice that $b-a$ fewer run components are used in the construction of $\phi(P)$ for $P\in\DD_{m,\al}^{(r+b,r)}$. If one use $\signature(P)$ for building the remaining steps of $\phi(P)$, the path has to go below diagonal without making up the 2 missing vertical components. Thus the construction ``rise up" the path by inserting signature 1 immediately before the path going below diagonal. The ``rising up" action will take place exactly $b-a$ times since before the $(b-a)$th rising up, there are still too few vertical components and the path is going below the diagonal following $\signature(P)$; after the $(b-a)$th rising up, the remaining signature are exactly that of $P$, saying no more return until the end of the construction. In the example in \fref{phieg}, the two ``rising up"s are marked green with their signature circled, and one can observe that the path after the second ``rising up", the remaining steps of $P'$ are the same as $P$.
		
		To show the bijectivity of $\phi$, we need to show that $\phi$ is surjective and injective. The injectivity is obvious since different {\LiLan plain} signatures will result in different images. The surjectivity follows from the fact that $\phi$ is invertible. 

        {\LiLan For any $X\in D^{(r+a,r+b-a)}_{m,\alpha}$ with a plain signature $\psig(X)=s_1\cdots s_{m-r-b-1}$, we construct an $m\times n$ lattice path $P$ with composition type $\al$ and signature $\sig(P)=1^{r+b}0s_1\cdots s_{m-r-b-1}$. We shall prove that $P$ is the preimage of $X$ under the map $\phi$, i.e.,\ $P\in D^{(r+b,r)}_{m,\alpha}$ with the same plain signature.

        A lattice path on the $m\times n$ grid is uniquely determined by its composition type $\al$ and its full signature.
        This path is a Dyck path if and only if it never crosses the main diagonal; algebraically, this means that for every $i=1,\ldots,\ell$ the $i$th $1$ in the signature must be followed by at least
        $
        \lceil\frac{m}{n}(\al_{i}+\cdots+\al_{\ell})\rceil -\ell+i-1
        $
        zeros.
        Moreover, the $i$th vertical segment is a return component if and only if the $i$th $1$ is followed by exactly
        $
        \lceil\frac{m}{n}(\al_{i}+\cdots+\al_{\ell})\rceil -\ell+i-1
        $
        zeros.

%         The plain signature $\psig(X)$ is a subsequence of $\sig(X)$, and
% \[
% \sig(X)=1^{r+a}0s'_1\dots s'_{m-r-a-1},
% \]
% where $s'_1\dots s'_{m-r-a-1}$ is a shuffle of $s_1\dots s_{m-r-b-1}$ and $1^{b-a}$; the latter $b-a$ ones correspond exactly to return components.
We may write $\sig(X)=1^{r+a}0s'_1\dots s'_{m-r-a-1}$, where $s'_1\dots s'_{m-r-a-1}$ is a shuffle of $s_1\dots s_{m-r-b-1}$ and $1^{b-a}$, the latter $b-a$ ones correspond exactly to the east steps in return components.
Thus $\sig(P)$ is obtained by moving these $b-a$ ones to the far left, so for every $i$ the $i$th $1$ in $\sig(P)$ is weakly to the left of the $i$th $1$ in $\sig(X)$; consequently, the resulting path $P$ is guaranteed to be a Dyck path.

We perform the move from left to right, one $1$ at a time, ensuring that no new return is created beyond column $r+b$.
Let $s'_i$ be the leftmost $1$ in $s'_1\dots s'_{m-r-a-1}$ that corresponds to a return component, and suppose it is the $j$th $1$ in $\sig(X)$.
Since $s'_i$ is leftmost, its immediate left neighbour in $\sig(X)$ must be a $0$; otherwise the $1$ one position further left would already mark a return component.
Moving $s'_i$ to the far left yields the signature
$$
1^{r+a+1}0s'_1\dots\widehat{s'_i}\dots s'_{m-r-a-1}.
$$
Here the $j$th $1$ has strictly more zeros to its right than the $j$th $1$ in $\sig(X)$, so the corresponding vertical segment is no longer a return component.
For every $h$ ($r+a<h<j$) the $h$th $1$ also has weakly more zeros to its right, hence remains non-return (as $s'_i$ was the leftmost return-carrying $1$); the part after $s'_{i+1}$ is unchanged.
Thus, this single move removes the return corresponding to $s'_i$ and creates no new returns.

Repeating the procedure for the remaining leftmost return-carrying ones produces $\sig(P)$; the resulting path $P$ lies in $\DD_{m,\alpha}$, has run value $r+b$, and by induction its returns are exactly the $r$ contributions from the run components. $\psig(P)=\psig(X)=s_1\dots s_{m-r-b-1}$ also follows immediately.}\end{proof}

	The involution $\Phi$ in \tref{run} can be constructed from
	the bijections mentioned in the proof of \tref{run2},
	by pairing paths in $\DD_{m,\al}^{(k,\ell)}$ and $\DD_{m,\al}^{(\ell,k)}$ that have the same {\LiLan plain} signature. Thus, \tref{run} follows immediately. We note that this also generalizes and extends Li-Lin's result \cite[Theorem 2.3]{Li} to the rational case.
	We have the following corollary.
	\begin{cor}\label{symm:equ1}
		Given positive integers $m>n$ and a composition $\al\vDash n$, we have
		\begin{eqnarray*}
			\sum_{P\in \DD_{m,\al}}p^{\run(P)}q^{\ret(P)}&=&\sum_{P\in \DD_{m,\al}}p^{\ret(P)}q^{\run(P)},\\
			\sum_{P\in \DD_{m,n}}p^{\run(P)}q^{\ret(P)}&=&\sum_{P\in \DD_{m,n}}p^{\ret(P)}q^{\run(P)}.
		\end{eqnarray*}
	\end{cor}
	
	Recall the alternative formulas for the rr statistic given in \eqref{rr0} or \eqref{rr}. When $\lfloor \frac{m}{n} \rfloor\geq2$, $\rr(\al,m)=1$ for any composition $\al\vDash n$. As a consequence, we have the following.
	\begin{cor}\label{corodet}
		Given positive integers $m,n,k,\ell$ such that $\lfloor \frac{m}{n} \rfloor\geq2$, we have $|\DD_{m, n}^{(k,\ell)}|=|\DD_{m, n}^{(k+\ell-1,1)}|$.
	\end{cor}
	This result is valid for all $k$-Dyck paths for $k\geq 2$, though not valid for classical cases when $m=n$. {If one takes $|\DD_{m, n}^{(k,\ell)}|$ as the $(k,\ell)$-entry of an $n\times n$ matrix, it is seen to be a Hankel matrix (i.e., with constant skew-diagonals).} For example, we have
	\begin{center}
		$\{|\DD_{12, 4}^{(k,\ell)}|\}=\left[
		\begin{matrix}
			52&30&8&1\\
			30&8&1&0\\
			8&1&0&0\\
			1&0&0&0\\
		\end{matrix}\right],\qquad \{|\DD_{13, 5}^{(k,\ell)}|\}=\left[
		\begin{matrix}
			106&106&39&9&1\\
			106&39&9&1&0\\
			39&9&1&0&0\\
			9&1&0&0&0\\
			1&0&0&0&0\\
		\end{matrix}\right].$
	\end{center}
	
	In the general case when positive integers $k,\ell,m,n$ with $m>n$ and composition $\al\vDash n$ are given, let $r=\rr(\al,m)$, then the cardinality $|\DD_{m,\al}^{k,\ell}|=|\DD_{m,\al}^{k+\ell-r,r}|$ is the number of lattice paths from $(k+\ell-r+1, \al_1+\cdots+\al_{k+\ell-r})$ to $(m,n)$ with no north steps bordering the diagonal cells of $m\times n$ lattice{, \LiLan i.e.,\ the cardinality of the set $\mathcal{T}(m,\alpha,k+\ell-r)$}. This observation leads to the following alternative characterization.
	{\LiLan\begin{cor}
    Given positive integers $m,n,\ell$, and $r=\rr(\al, m)$. We have \[|\DD_{m,\al}^{k,\ell}|=|\DD_{m,\al}^{k+\ell-r,r}|=|\mathcal{T}(m,\alpha,k+\ell-r)|.
    \]
    \end{cor}}
	\begin{proof}
	Given a path $P\in\DD_{m,\al}^{k+\ell-r,r}$, recall the decomposition $P=M_1M_2\cdots M_{k+\ell-r}ET$, where $M_i=N^{\al_i}E$ for $1\le i\le k+\ell-r$. Now the path $T$ is indeed from $(k+\ell-r+1, \al_1+\cdots+\al_{k+\ell-r})$ to $(m,n)$ and it must not border any diagonal cells to ensure that no more return positions are created. Conversely, taking any such partial lattice path as the suffix, one uses the values of $k+\ell-r$ and $\alpha$ (its first $k+\ell-r$ parts to be precise) to construct the prefix and recover the original path $P\in\DD_{m,\al}^{k+\ell-r,r}$.
    % must ends with an $E$ step, we can write $T=\tilde{T}E$. Now shifting the entire partial path $\tilde{T}$ one unit to the right yields the desired $\tilde{P}$. The whole process is seen to be invertible. 
    \end{proof}

{\bf \noindent Remark 1.}
The number of paths satisfying $\run(P)=k,\ \ret(P)=\ell$ and $\comp(P)=\al$ equals the number of paths satisfying $\run(P)=k+\ell-r,\ \comp(P)=\al$ with no returns after the dashed $E$ step. Therefore, the set of valid {\LiLan plain signatures for $\DD_{m,\al}^{k,\ell}$ corresponds to the set of signatures of paths in $\mathcal{T}(m,\alpha,k+\ell-r)$. }

As an example, for Dyck paths with $\run(P)=1,\ \ret(P)=2$ and $\comp(P)=2111$ in $11\times 5$ lattice, we know that they are in bijection with paths with $\run(P)=2,\ \ret(P)=1$ and $\comp(P)=2111$. Thus, their valid {\LiLan plain} signatures correspond to the signatures of paths weakly bounded by the two gray paths shown in \fref{valid_sig} from lattice point $(3,3)$ to $(11,5)$. It is a pure enumeration problem to compute the number of valid {\LiLan plain} signatures, and we shall leave this problem open to interested readers.

\begin{figure}[ht!]
		\centering
		\begin{tikzpicture}[scale=0.39]
			\mngrid{0,0}{11}{5}
			\colorpath{0,0}{2,-1,1,-1}{black}
            \colorpath{2,3}{-1}{red,dotted}
            \colorpath{3,3}{1,-1,1,-1,-1,-1,-1,-1,-1,-1}{black!30}
            \colorpath{3,3}{-1,-1,1,-1,-1,1,-1,-1,-1,-1}{black!30}
		\end{tikzpicture}
		\caption{$P\in \DD_{11, 2111}^{2,\ell}$}
		\label{fig:valid_sig}
	\end{figure}

	\subsection{\texorpdfstring{The symmetry of $(\runb,\ret)$ distribution on $\DD_{m,n}$}{The symmetry of (run,ret) distribution on Dmn}}
	For the ratio-run statistic $\runb$, we have a similar result.
	\begin{thm}\label{theorem:runb1}
		Given $m>n$, there is a composition type preserving involution $\Psi$ on the set $\DD_{m,n}$, sending statistics $(\runb, \ret)$ to $(\ret, \runb)$.
	\end{thm}
	
	Let 
	$$\tDD_{m,n}^{(a,b)}:=\{P\in \DD_{m,n}:\runb(P)=a,\ \ret(P)=b\},$$
	and
	$$\tDD_{m,\al}^{(a,b)}:=\{P\in \DD_{m,\al}: \runb(P)=a,\ \ret(P)=b\}.$$
	  Given a path $P\in\tDD_{m,\al}^{(a,b)}$, let $n=|\al|$ and $k = \lfloor \frac{m}{n} \rfloor$. Notice that each $\runb$ statistic of $P$ corresponds 
to a block of $k$ east steps, specifically those from column $ik+1$ to $(i+1)k$, $(0\leq i\leq a-1)$.
    %Given a path $P\in\tDD_{m,\al}^{(k,\ell)}$, we set $n=|\al|$ and $k=\lfloor \frac{m}{n}\rfloor$. Notice that each $\runb$ statistic of $P$ corresponds to $k$ east steps with signature that begins with 1. 
    We shall define signature statistics not only for steps not contributing to $\runb$ and return (like the normal run statistic), but also for steps contributing to $\runb$ and return. 
    We first collect steps contributing to either {$\runb$} or $\ret$ or both as $P_1$, i.e., we scan the path from the beginning, whenever a {$\runb$} or a return appear, we collect $k$ steps and append to $P_1$; the path $P_2$ consists of the remaining steps except the dashed east step. For the left path $P$ in \fref{runbex}, we have $P_1=NENENNENENENENEE$ appearing in columns 1,2,3,4,13,14,21,22, and $P_2=EENEEEEENENEEEEEE$ appearing in columns 6,7,8,9,10,11,12,15,16,17,18,19,20,23. In fact, the path $P$ is a shuffle of $P_1$ and $P_2$. {For a sequence $X$ of steps, let $\sig(X)$ be the signature of its east steps.} 
	We define the \textit{ratio-signatures} of $P$ to be
	\begin{eqnarray*}
		\signa(P) &=& \sig(P_1),\\
		\signb(P) &=& \sig(P_2).
	\end{eqnarray*} 
	{For instance, the left path $P$ in \fref{runbex} has two ratio-runs (in lattice columns 1 and 2, 3 and 4) and three returns (in lattice columns 1, 13, 21).
    We collect the steps in columns 1,2,3,4,13,14,21,22 as $P_1$, the step in column 5 is a forced east step that is not counted in signature, and the remaining steps form path $P_2$.
    We have $\sig(P_1)=11111110$ (black) and $\sig(P_2)=00100001100000$ (blue)}. %{ Therefore, for the $\sig(P)=(\underline{{\color{green}1,1},{\color{red}1,1},0,0,0,1,0,0,0,0},\underline{1,1,1,1,0,0,0,0},\underline{1,0,0})$, we note that
%the two monochromatic pairs {\color{green}1,1} and {\color{red}1,1} each generate one run (totaling 2 runs), since a consecutive subsequence of ones at the beginning with length $k=2$ is counted as a run under our definition. And each of the three underlined segments contributes one ret, resulting in 3 rets in total. There is a similar discussion about $\sig(P')=(\underline{{\color{green}1,1},0,0},\underline{1,1,0,1,0,0,0,0},\underline{1,1,1,1,0,0,0,0},\underline{1,0,0})$.}
    
	We let $\trr(P)$ denote the largest integer $i$ such that the first $i$ vertical components all contribute to both $\runb$ and ret statistics, where $\trr(P)-1$ is the number of ratio-runs bordering the diagonal cells.
	
	Let $P$ be a path in $\DD_{m,n}$, the statistic $\trr$ depends on the composition type $\comp(P)=\al\vDash n$, the number $k$ and $\signa(P)=s$. We shall also write $\trr(\al,m,s)$ for such $P\in\DD_{m,n}$.
	We shall omit the formula for $\trr$, but only state that $\trr(P)-1$ is the maximal initial consecutive occurrence of path pattern $NE^k$ such that the south ends of the $N$ steps are return positions.
	If we let
	$$\DD_{m,\al, s}:=\{P\in \DD_{m,\al}: \signa(P)=s\}$$
	and
	$$\tDD_{m,\al, s}^{(a,b)}:=\{P\in \DD_{m,\al}: \runb(P)=a,\ \ret(P)=b,\ \signa(P)=s\},$$
	then $\trr(P)$ is the smallest number that $\runb(P)$ or $\ret(P)$ can reach for $P\in\DD_{m,\al, s}$.

	The proof of \tref{runb1} is analogous to \tref{run}, where we need the following stronger result.
	\begin{thm}\label{theorem:runb2}
		Let $m>n$ be positive integers, $\al = (\al_1,\ldots,\al_\ell)$ be a composition of $n$ and $s$ be a 01-sequence. If $r=\trr(\al,m,s)$, then
		\begin{equation}\label{runbthm}
			\left|\tDD_{m,\al,s}^{(r+a,r+b-a)}\right| = \left|\tDD_{m,\al,s}^{(r+b,r)}\right|
		\end{equation}
		for any non-negative integers $a\leq b$.
	\end{thm}
    
	\tref{runb1} follows from \tref{runb2}, and the proof of \tref{runb2} is completely analogous to that of \tref{run2}, {the key distinction from the proof of \tref{run2} lies in the fact that each $M_i$ in the decomposition of $P$ contains $k$ east steps}. Thus, we can build a series of bijective maps among sets of the form $\tDD_{m,\al,s}^{(x,y)}$. The idea of the proof is that given a path with certain signatures from the first set, there is a unique path from the second set having the same signatures.
    We shall omit the details and give an example. In \fref{runbex}, we have a path $P\in\tDD_{23,(1,1,2,1,1,1,1,1,1,1),11111110}^{(2,3)}$, and we construct $P'\in\DD_{23,(1,1,2,1,1,1,1,1,1,1),11111110}^{(1,4)}$ with the same signatures. In the path $P'$, the only ratio-run appears in columns 1 and 2, and the 4 returns appear in columns 1, 5, 13 and 21. The part of path that corresponds to $\signa(P')$ is marked black in the figure.
	
	\begin{figure}[ht!]
		\centering
		\begin{tikzpicture}[scale=0.32]
			\mngrid{0,0}{23}{11}
			\colorpath{0,0}{1,-1,1,-1,2,-1,1,-1}{}
			\signseq{0,0}{2/1,3/1,5/1,6/1}
			\colorpath{4,5}{-1}{red,dotted}
			\colorpath{5,5}{-2,1,-5}{cyan}
			\signseqb{5,5}{1/0,1/0,2/1,2/0,2/0,2/0,2/0}
			\colorpath{12,6}{1,-1,1,-1}{}
			\signseq{12,6}{2/1,3/1}
			\colorpath{14,8}{1,-1,1,-5}{cyan}
			\signseqb{14,8}{2/1,3/1,3/0,3/0,3/0,3/0}
			\colorpath{20,10}{1,-2}{}
			\signseq{20,10}{2/1,2/0}
			\colorpath{22,11}{-1}{cyan}
			\signseqb{22,11}{1/0}

			\mngrid{27,0}{23}{11}
			\colorpath{27+0,0}{1,-1,1,-1}{}
			\signseq{27+0,0}{2/1,3/1}
			\colorpath{27+2,2}{-1}{red,dotted}
			\colorpath{27+3,2}{-1}{cyan}
			\signseqb{27+3,2}{1/0}
			\colorpath{27+4,2}{2,-1,1,-1}{}
			\signseq{27+4,2}{3/1,4/1}
			\colorpath{27+6,5}{-1,1,-5}{cyan}
			\signseqb{27+6,5}{1/0,2/1,2/0,2/0,2/0,2/0}
			\colorpath{27+12,6}{1,-1,1,-1}{}
			\signseq{27+12,6}{2/1,3/1}
			\colorpath{27+14,8}{1,-1,1,-5}{cyan}
			\signseqb{27+14,8}{2/1,3/1,3/0,3/0,3/0,3/0}
			\colorpath{27+20,10}{1,-2}{}
			\signseq{27+20,10}{2/1,2/0}
			\colorpath{27+22,11}{-1}{cyan}
			\signseqb{27+22,11}{1/0}
			
			\node at (25,5.5) {$\longrightarrow$};
		\end{tikzpicture}
		\caption{$P\in\tDD_{23,(1,1,2,1,1,1,1,1,1,1),11111110}^{(2,3)}$ and $P'\in\DD_{23,(1,1,2,1,1,1,1,1,1,1),11111110}^{(1,4)}$.}
		\label{fig:runbex}
	\end{figure}
	
	We also have the following corollary.
	\begin{cor}\label{symm:equ2}
		Given positive integers $m>n$, a composition $\al\vDash n$ and a 01-sequence $s$, we have
		\begin{eqnarray*}
			\sum_{P\in \DD_{m,\al,s}}p^{\runb(P)}q^{\ret(P)}&=&\sum_{P\in \DD_{m,\al,s}}p^{\ret(P)}q^{\runb(P)},\\
			\sum_{P\in \DD_{m,\al}}p^{\runb(P)}q^{\ret(P)}&=&\sum_{P\in \DD_{m,\al}}p^{\ret(P)}q^{\runb(P)},\\
			\sum_{P\in \DD_{m,n}}p^{\runb(P)}q^{\ret(P)}&=&\sum_{P\in \DD_{m,n}}p^{\ret(P)}q^{\runb(P)}.
		\end{eqnarray*}
	\end{cor}
	
	Setting $p=1$ in both Corollaries~\ref{symm:equ1} and \ref{symm:equ2}, it follows that the statistics $\ret, \run$ and $\runb$ are equidistributed on $\DD_{m,n}$ and $\DD_{m,\al}$, i.e.
	\begin{cor}
		Given positive integers $m>n$ and a composition $\al\vDash n$, we have
		\begin{eqnarray*}
			\sum_{P\in \DD_{m,\al}}q^{\ret(P)}&=&\sum_{P\in \DD_{m,\al}}q^{\run(P)}\ \ =\ \ \sum_{P\in \DD_{m,\al}}q^{\runb(P)},\\
			\sum_{P\in \DD_{m,n}}q^{\ret(P)}&=&\sum_{P\in \DD_{m,n}}q^{\run(P)}\ \ =\ \ \sum_{P\in \DD_{m,n}}q^{\runb(P)}.
		\end{eqnarray*}
	\end{cor}

	\section{Generating function results}\label{gen}
	In this section, we shall give generating function results for classical Dyck paths and
	$k$-Dyck paths involving the statistics $\run$ and $\ret$. {The generating functions not only prove the symmetries of statistics in some cases, but also imply many other statistical or analytical properties of statistics.}
	\subsection{Classical Dyck path case}
        {In the classical case, Li and Lin \cite{Li} gave a bijective proof that the statistics (run, ret) have symmetric distribution on classical Dyck paths. We shall give a generating function proof of this symmetry.} Given a path $P\in\DD_{n}$, let $\lvl(P)$ be the number of vertical components of $P$, and let
	$C(t,p,q,r)$ be the generating function:
	$$
	C(t,p,q,r):=\sum_{n\geq 0} t^n\sum_{P\in\DD_{n}} p^{\run(P)} q^{\ret(P)} r^{\lvl(P)}.
	$$
	$C(t,p,q,r)$ is a generalization of the generating function of Catalan numbers. We compute $C(t,p,q,r)$ as follows.
	
	Given a Dyck path $D$,
	\begin{enumerate}[label=(\roman*)]
		\item $D=\emptyset$, contributing 1 to $C(t,p,q,r)$.
		\item $D=NED'$, where $D'$ is a Dyck path, contributing $tpqr\cdot C(t,p,q,r)$ to $C(t,p,q,r)$.
		\item $D=ND_1ED_2$, where $D_1$, $D_2$ are Dyck paths and $D_1\neq\emptyset$, contributing $tq(C(t,p,1,r)-1)\cdot C(t,1,q,r)$ to $C(t,p,q,r)$.
	\end{enumerate}
	Thus, we have
	$$
	C(t,p,q,r)=1+tpqr\cdot C(t,p,q,r)+tq\cdot (C(t,p,1,r)-1)\cdot C(t,1,q,r).
	$$
	We let $C_p=C(t,p,1,r)$, $C_q=C(t,1,q,r)$ and $C_r=C(t,1,1,r)$ be abbreviations of the specializations of $C(t,p,q,r)$, then following a similar analysis, we have
	\begin{eqnarray*}
		C_r&=&1+tr\cdot C_r+t\cdot (C_r-1)\cdot C_r,\\
		C_p&=&1+tpr\cdot C_p+t\cdot (C_p-1)\cdot C_r,\\
		C_q&=&1+tqr\cdot C_q+tq\cdot (C_r-1)\cdot C_q.
	\end{eqnarray*}
	First solving $C_r$, then $C_p$ and $C_q$, finally $C(t,p,q,r)$, we have 
	\begin{eqnarray*}
		C_r&=&\frac{1-r t+t-\sqrt{(-r t+t+1)^2-4 t}}{2 t},\\
		C_p&=&\frac{2}{2-prt+pt-p+ p \sqrt{(r-1)^2 t^2-2 (r+1) t+1}},\\
		C_q&=&\frac{2}{2-qrt+qt-q+ q \sqrt{(r-1)^2 t^2-2 (r+1) t+1}},
	\end{eqnarray*}
	and
	\begin{eqnarray*}
		C(t,p,q,r)=\frac{2 p q (r^2t^2-r t^2-t+1)+2(p+q)(t- r t-1)+4+2(p + q - p q - p q r t) \sqrt{\Delta} }{(1-p q r t)\left(2-prt+pt-p+p \sqrt{\Delta}\right) \left(2-qrt+qt-q+q \sqrt{\Delta}\right) },
	\end{eqnarray*}
	where $\Delta:=(r-1)^2 t^2-2(r+1) t+1$. Clearly, $C(t,p,q,r)$ is symmetric in $p$ and $q$. {Let $r=1$, we can state the following theorem.
    \begin{thm}
        Let $C(t,p,q):=\sum_{n\geq 0} t^n\sum_{P\in\DD_{n}} p^{\run(P)} q^{\ret(P)}$ be the generating function of run and return. Then
        $$
        C(t,p,q)=
		\frac{2 p q (1-t)-2(p+q)+4+2(p + q - p q - p q t) \sqrt{ -4t+1} }{(1-p q  t)\left(2-p+p \sqrt{ -4t+1}\right) \left(2-q+q \sqrt{ -4 t+1}\right) }.
        $$
    \end{thm}
    }

    One can compute by Mathematica to obtain the Taylor series of the generating function $C(t,p,q,r)$ that
    \begin{eqnarray*}
        &&C(t,p,q,r)=1+p q r t+t^2 \left(p^2 q^2 r^2+p q r\right)
        +t^3 \left(p^3 q^3 r^3+r^2
        \left(p^2 q^2+p^2 q+p q^2\right)+p q r\right)
        \\&&+t^4
        \left(p^4 q^4 r^4+r^2 \left(p^2 q^2+2 p^2 q+2 p q^2+p
        q\right)+r^3 \left(p^3 q^3+p^3 q^2+p^3 q+p^2 q^3+p^2
        q^2+p q^3\right)+p q r\right)
        \\&&+t^5 \left(p^5 q^5 r^5
        +r^4 \left(p^4 q^4+p^4 q^3+p^4
        q^2+p^4 q+p^3 q^4+p^3 q^3+p^3 q^2+p^2 q^4+p^2 q^3+p
        q^4\right)
        \right.
        \\&&\left.\qquad +r^3 \left(p^3
        q^3+2 p^3 q^2+3 p^3 q+2 p^2 q^3+4 p^2 q^2+2 p^2 q+3 p
        q^3+2 p q^2+p q\right)\right.
        \\&&\left.\qquad +r^2
        \left(p^2 q^2+3 p^2 q+3 p q^2+3 p q\right) +p q r\right)+\cdots.
    \end{eqnarray*}
Some coincidence may be found from the On-Line Encyclopedia of Integer Sequences (OEIS) \cite{OEIS} if one does some manipulations, such as specializations, substitutions, or looking at coefficients of some special terms of the series. 
    
In addition to the series, the generating function itself gives rich analytical information of the statistics. For example, 
a partial derivative $\partial C(t,p,q) /\partial p$ and specializing at $p=0$ gives
\begin{eqnarray*}
F(t,q):=\frac{C(t,p,q)}{\partial p}\bigg|_{p=0} = \frac{(q^2 - q)t\sqrt{1 - 4t} - (q^2 - 3q)t}{q\sqrt{1 - 4t} - q + 2},
\end{eqnarray*}
which is the generating function of Dyck paths whose run is 1, and the power of $q$ still tracks the return of paths.
The Taylor series of $F(t,q)$ is
\begin{eqnarray*}
F(t,q) = qt + qt^2 + (q^2+q)t^3 + (q^3+2q^2+2q)t^4 + (q^4+3q^3+5q^2+5q)t^5 + O(t^6). 
\end{eqnarray*}
One may look up OEIS \cite{OEIS} to find this as sequence A009766, which is the Catalan's triangle $T(n,k)$.

Furthermore, taking the partial derivative of $F(t,q)$ with respect to $q$ and evaluating at $q=1$ yields 
\begin{eqnarray*}
G(t) :=\frac{F(t,q)}{\partial q}\bigg|_{q=1}  = \frac{2(t^2 - t)}{2t - \sqrt{1 - 4t} - 1}.
\end{eqnarray*}
This is the generating function enumerating the total sum of returns of paths with exactly 1 run. Its Taylor series is:
\begin{eqnarray*}
F(t) = t + t^2 + 3t^3 + 9t^4 + 28t^5 + 90t^6 + 297 t^7 + 1001 t^8 + O(t^9),
\end{eqnarray*}
which appears in OEIS \cite{OEIS} as sequence A071724, enumerating the number of standard Young tableaux of shape $(n+1,n-1)$. We leave it to the interested reader to provide bijective justifications of these connections with the two OEIS entries.

	\subsection{\texorpdfstring{$k$-Dyck path case}{k-Dyck path case}}
	Given a path $P\in\DD_{kn,n}$, let $\lvl(P)$ be the number of vertical components of $P$, and let
	$C^{(k)}(t,p,q,r)$ be the generating function:
	$$
	C^{(k)}=C^{(k)}(t,p,q,r):=\sum_{n\geq 0} t^n\sum_{P\in\DD_{kn,n}} p^{\run(P)} q^{\ret(P)} r^{\lvl(P)}.
	$$
	%We first compute $C(t,p,q,r)$ as follows.
	
	We decompose a Dyck path $D$ in a similar manner as the previous subsection.
	\begin{enumerate}[label=(\roman*)]
		\item $D=\emptyset$, contributing 1 to $C^{(k)}$.
		\item $D=NED_2...D_{k+1}$, where $D_i(i>1)$ is a Dyck path, contributing $tpqr\cdot C^{(k)}(t,p,1,r)\cdot(C^{(k)}(t,1,1,r))^{k-2}\cdot C^{(k)}(t,1,q,r)$ to $C^{(k)}$.
		\item $D=ND_1ED_2...ED_{k+1}$, where $D_i(i>0)$ are Dyck paths and $D_1\neq\emptyset$, contributing $tq(C^{(k)}(t,p,1,r)-1)\cdot (C^{(k)}(t,1,1,r))^{k-1}\cdot C^{(k)}(t,1,q,r)$ to $C^{(k)}$.
	\end{enumerate}
	We let $C^{(k)}_p=C^{(k)}(t,p,1,r)$, $C^{(k)}_q=C^{(k)}(t,1,q,r)$ and $C^{(k)}_r=C^{(k)}(t,1,1,r)$ be abbreviations of the specializations of $C^{(k)}$, we have
	\begin{align*}
		C^{(k)}&=1+tq C^{(k)}_q (C^{(k)}_p C^{(k)}_r + prC^{(k)}_p - C^{(k)}_r) (C^{(k)}_r)^{k-2},\\
		C^{(k)}_p&=1+t(C_p^{(k)} C_r^{(k)}+prC_p^{(k)}-C_r^{(k)})(C^{(k)}_r)^{k-1},\\
		C^{(k)}_q&=1+tq C_q^{(k)} (C^{(k)}_r+r-1)(C^{(k)}_r)^{k-1},\\
		C^{(k)}_r&=1+t(C^{(k)}_r+r-1)(C^{(k)}_r)^k.
	\end{align*}
	Thus,
	\begin{equation*}
		C^{(k)}=\frac{C^{(k)}_r\left((p-1)(q-1)C^{(k)}_r + (p+q)-pq\right)}{\left((p-1)C^{(k)}_r - p\right)\left((q-1)C^{(k)}_r-q\right)}.
	\end{equation*}
	From this expression, one can see that $C^{(k)}$ is symmetric in $(p,q)$.

For the $k$-Dyck path case, one can still obtain Taylor series from the recursive relations of the generating functions. On the generating function manipulation side, one may notice that the generating function $C^{(k)}_r$ has a nice recursive formula, which allows one to utilize Lagrange inversion to compute coefficients; and the more comprehensive generating function, $C^{(k)}$, can be written in terms of $C^{(k)}_r$. One may try to use analytical or algebraic tools to explore the generating functions and properties of these statistics.

	\section{Concluding remarks}\label{open}
	In this paper, we show that two pairs of statistics have symmetric joint distributions on the set of rational Dyck paths. There are still a large number of {unsolved enumeration problems along this line of inquiry}, such as {finding} an explicit formula for the number of rational Dyck paths with $k$ runs (of two types) and $\ell$ returns. It is also possible to generalize our results to other sets of combinatorial objects, such as lattice paths and partial Dyck paths. 
    
   Our generalized results in the rational case have the restriction that $m \geq n$, when the statistics run and return have symmetric distribution. On the other hand when $m<n$, run and return do not have symmetric distribution under the current definition of statistics. For instance, the Dyck path $P = NNNENENNEE$ in a $4 \times 6$ lattice satisfies $\run(P) = 3$  and $\ret(P) = 2$, while no path in the same lattice has $\run(P) = 2$ and $\ret(P) = 3$. Notice that our bijections require the condition $m\ge n$ to ensure that every such path terminates at a proper step. In order to obtain symmetric distribution results for rational cases when $m<n$, people have to modify the run and return statistics. 
    So we have the open problem that for any integers $m$ and $n$, modify the definition of run and return such that they have symmetric joint distributions.

	Let Inor be a statistic of Dyck paths that counts the number of initial consecutive north steps. {In \cite{Li}, Li and Lin proved that the joint distribution of (run, Inor) and (ret, Inor) are symmetric for classical Dyck paths, where the latter was first established in \cite{Deu1}. 
    Elizalde and Rubey directly establishes the symmetric joint distribution of (ret, Inor) for rational Dyck paths (Theorem 2.1 in \cite{Eli}). 
 %   Also, Theorem 2.1 in \cite{Eli}).  directly establishes the symmetric joint distribution of (ret, Inor) for rational Dyck paths - where the maximal path $N^m E^n$ tightly bounds the diagonal cells.
    We are able to} extend both results to the rational case when $m\geq n$, where the bijection we give for the pair (ret, Inor) is different from the bijection of Elizalde and Rubey. Our relevant work is in preparation.
	
	It is worth noting that both of our bijections preserve composition type, which is stronger than preserving partition type. Another well-studied combinatorial object is parking function that we mentioned in Section~1. Rational parking functions are labeled rational Dyck paths where labels in each column are increasing from bottom to top. Let $\mathcal{PF}_{m,n}$ be the set of rational $(m,n)$-parking functions, and let the run and return statistic of a parking function be that of the underlying Dyck path. As a consequence, the two generating functions, 
	\begin{eqnarray*}
		\sum_{P\in\mathcal{PF}_{m,n}} p^{\run(P)}q^{\ret(P)}X^P &=& \sum_{P \in \DD_{m,n}}p^{\run(P)}q^{\ret(P)}e_{\lambda(P)},\\
		\sum_{P\in\mathcal{PF}_{m,n}} p^{\runb(P)}q^{\ret(P)}X^P &=& \sum_{P \in \DD_{m,n}}p^{\runb(P)}q^{\ret(P)}e_{\lambda(P)}
	\end{eqnarray*}
	are symmetric functions with coefficients symmetric in $(p,q)$, where $e_\la$ is the elementary symmetric function. Let $\nabla$ be the eigen-operator of the Macdonald polynomials that $\nabla \tilde{H}_\mu = T_\mu \tilde{H}_\mu.$
	When $p=q=1$ and $m=n$, the generating functions above become $\nabla e_n|_{q=t=1}$, which is a specialization of the shuffle theorem of Carlsson and Mellit \cite{Car}. This may also hold potential interest for the algebraic combinatorics community.

	\section*{Acknowledgement}
	
	We would like to express our sincere gratitude to the reviewer for insightful comments and constructive suggestions, which have significantly strengthened the manuscript. We thank Martin Rubey for helpful comments, and Yang Li and Zhicong Lin for insightful discussions that improved this work. This work is supported by the Fundamental Research Funds for the Central Universities. Shishuo Fu is supported by the National Natural Science Foundation of China (12171059 and 12371336) and the Mathematical Research Center of Chongqing University.
	Dun Qiu is supported in part by the National Natural Science Foundation of China (12271023), and the Natural Science Foundation of Tianjin (24JCZDJC01390).

\end{document}